\documentclass [10pt,reqno]{amsart}
\usepackage{ucs}
\usepackage{amsmath, amssymb, amscd}
\usepackage{pb-diagram}
\usepackage[latin1]{inputenc}
\usepackage{amsmath}
\usepackage{amssymb}

\newtheorem{theorem}{Theorem}[section]
\newtheorem{definition}[theorem]{Definition}
\newtheorem{lemma}[theorem]{Lemma}
\newtheorem{prop}[theorem]{Proposition}

\newtheorem*{remark}{Remark}

\newcommand{\abs}[1]{\lvert#1\rvert}
\newcommand{\norm}[1]{\| #1 \|}
\newcommand{\klammer}[1]{\left( #1 \right)}

\renewcommand{\epsilon}{\varepsilon}

\DeclareMathAlphabet{\mathpzc}{OT1}{pzc}{m}{it}
\usepackage{mathrsfs}

\newcommand{\Z}{\mathbb{Z}}
\newcommand{\C}{\mathbb{C}}

\newcommand{\R}{\mathbb{R}}

\renewcommand{\qed}{$\hfill \square$ \bigskip \\}
\renewcommand{\phi}{\varphi}

\newcommand{\ad}{\text{ad}}
\newcommand{\Ad}{\text{Ad}}
\newcommand{\rk}{\text{rk}}

\newcommand{\tensor}{\otimes}
\newcommand{\End}{\text{End}}
\newcommand{\Hom}{\text{Hom}}
\newcommand{\Aut}{\text{Aut}}
\renewcommand{\det}{\text{det}}
\newcommand{\id}{\text{id}}

\newcommand{\Atheta}{\mathscr{A}_{\theta}(P_E)}

\newcommand{\su}{\mathfrak{su}}
\renewcommand{\sl}{\mathfrak{sl}}
\newcommand{\gl}{\mathfrak{gl}}
\renewcommand{\u}{\mathfrak{u}}

\renewcommand{\bar}{\overline}
\newcommand{\G}{\mathscr{G}^0}

\newcommand{\tr}{\text{tr}}
\newcommand{\data}{\mathfrak{s},E}
\newcommand{\dataf}{\mathfrak{s},F}
\newcommand{\conf}{\mathscr{C}}
\newcommand{\bonf}{\mathscr{B}}
\newcommand{\s}{\mathfrak{s}}

\newcommand{\hata}{{\hat{A}}}

\begin{document}

\thispagestyle{empty}

\title[Monopoles, higher rank instantons, and Seiberg-Witten invariants]{$PU(N)$ monopoles, higher rank instantons, and the monopole invariants}
 \author
[Raphael
Zentner]{Raphael Zentner} 

\begin {abstract} 
A famous conjecture in gauge theory mathematics, attributed to Witten, suggests that the polynomial invariants of Donaldson are expressible in terms of the Seiberg-Witten invariants if the underlying four-manifold is of simple type. Mathematicians have sought a proof of the conjecture by means of a `cobordism program' involving $PU(2)$ monopoles. 
A higher rank version of the Donaldson invariants was recently introduced by Kronheimer. Before being defined, the physicists Mari\~no and Moore had already suggested that there should be a generalisation of Witten's conjecture to this type of invariants. We adopt a generalisation of the cobordism program to the higher rank situation by studying $PU(N)$ monopoles. We analyse the differences to the $PU(2)$ situation, yielding evidence that a generalisation of Witten's conjecture should hold. 
\end {abstract}

\address {Fakult\"at f\"ur Mathematik \\ Universit\"at Bielefeld \\ 33501 Bielefeld\\
Germany}

\email{rzentner@math.uni-bielefeld.de}

\maketitle

\section*{Introduction}
The two important gauge-theoretical invariants of a smooth closed four-manifold are the polynomial invariants of Donaldson \cite{D}, derived from anti-selfdual $PU(2)$ connections in rank-2-bundles, and the Seiberg-Witten invariants, derived from the Seiberg-Witten equations \cite{W} which are associated to $Spin^c$ structures on the four-manifold. 
Kronheimer and Mrowka have proved a structure theorem for the
Donaldson invariants \cite{KM2}, showing that for 4-manifolds of `simple type'
the polynomial invariants are specified by certain algebraic-topological data, in particular the intersection form,
and a finite set of distinguished cohomology
classes in the group $H^2(X,\Z)$ each coming with some rational coefficient. 
With this at hand, Witten claims that the polynomial invariants are determined
by the Seiberg-Witten-invariants, with the basic classes being the first Chern classes
of the $Spin^c$-structures with non-trivial Seiberg-Witten-invariant, and with an explicit formula for
the rational coefficients \cite{W}. Witten derived this conjecture from correlation functions in quantum field theory and certain limiting behaviours with respect to a certain coupling constant. Mathematicians then tried to derive a proof of the conjecture by a certain cobordism obtained from $PU(2)$ monopoles. Heuristically the idea is as follows (ignoring the technical problems involved). There is a circle action on the moduli space of $PU(2)$ monopoles which comes from multiplying the spinor component by a complex number of norm one. The fixed point set of this action consists of the moduli space of $PU(2)$ instantons, and further a finite number of moduli spaces of classical $U(1)$ Seiberg-Witten monopoles. The circle acts freely on the complement of this fixed point locus, and the quotient yields a cobordism between a projective bundle over the moduli space of $PU(2)$ instantons, and projective bundles over the moduli spaces of $U(1)$ Seiberg-Witten monopoles. Furthermore, the canoncial cohomology class that yields the polynomial invariant of Donaldson extends to the cobordism. The evaluation of this extension on one side, yielding the Donaldson invariant, is therefore equal to the evaluation on the other sides, which should be expressions containing the Seiberg-Witten invariants of the moduli spaces in the fixed point locus. This program was started independently by Pidstrigach and Tyurin \cite{PT} and Okonek and Teleman \cite{OT2}, \cite{T}. It was carried on over years by Feehan and Leness \cite{FL1}, \cite{FL2}, \cite{FL3}, \cite{FL4}. It seems that they have now proved the full conjecture \cite{FL5}.

Recently Kronheimer introduced polynomial invariants associated to moduli spaces of anti-selfdual $PU(N)$ connections in Hermitian rank-N-bundles \cite{K}. These are generalisations of the Donaldson invariants, but the technical problems are much harder than in the classical situation. Before these invariants were even properly defined, the physicists Mari\~no and Moore \cite{MM} had conjectured that there should be a generalisation of Witten's conjecture for these invariants, implying in particular that they do not contain new differential-topological information. Their argument relies again on physics. Kronheimer computed explicitely his higher rank invariants for manifolds obtained through knot surgery on the K3-surface, confirming the conjecture for this class of examples. Our approach to the situation is to introduce $PU(N)$ monopoles in order to follow the cobordism strategy described above. For a given $Spin^c$ structure $\mathfrak{s}$ and a Hermitian rank-N-bundle $E$ on $X$, the configuration space will consist of sections $\Psi$ of the `twisted spinor bundle' $W^+ = S^+_\mathfrak{s} \tensor E$ and by unitary connections $\hata$ in $E$ with fixed induced connection $\theta$ in the determinant line bundle of $E$. The straightforward generalisation of the $PU(2)$ monopole equations then read:
\begin{equation*}
\begin{split}
 D_{\hat{A}} \Psi & =  0 \\
 \gamma(F_A^+) - \mu_{0,0}(\Psi) & =  0 \ .
\end{split}
\end{equation*}
Here $D_{\hata}$ is the associated Dirac-operator to $\hata$, the map $\gamma$ is derived from Clifford-multiplication, $F_A^+$ is the self-dual part of the curvature of the $PU(N)$ connection $A$ induced by $\hata$, and $\mu_{0,0}$ is a quadratic map in the spinor which is explicitely described below. The gauge group of the problem is that of special unitary automorphisms of $E$. Again, there is a circle action on the moduli space of these $PU(N)$ monopoles which is given by the formula $(z,[\Psi,\hata]) \mapsto [z^{1/N} \Psi,\hata]$. The moduli space of $PU(N)$ instantons is contained as the locus of monopoles with vanishing spinor, and the other fixed point loci are labelled by a finite number of isomorphism classes of proper subbundles $[F]$ of $E$. An equivalence class $[\Psi,\hata]$ belongs to the $[F]$-locus $M^{[F]}$ if for each $F \in [F]$ there is a representative $(\Psi,\hata)$ with the spinor $\Psi$ being a section of $S^+_\mathfrak{s} \tensor F$ and with the connection $\hata$ splitting according to $E = F \oplus F^\perp$. It turns out that if $X$ is simply connected the description of $M^{[F]}$ is particularly simple after fixing one such $F$: The content of Theorem \ref{S1fixedpointset_main_text} is that we have a `parametrisation' 
\begin{equation}\label{contributions}
 M^U_{\mathfrak{s},F} \times M^{asd}_{F^\perp} \to M^{[F]} \ .
\end{equation}
Here $M^U_{\mathfrak{s},F}$ is a moduli space of $U(n)$ monopoles \cite{Z_un} with $n= \rk(F)$ having possible values $1 \leq n < N$, and $M^{asd}_{F^\perp}$ is the moduli space of anti-selfdual $PU(N-n)$ connections in $F^\perp$. In the case $n=1$ the moduli space $M^U_{\mathfrak{s},F}$ is a classical $U(1)$ Seiberg-Witten moduli space. The map (\ref{contributions}) is surjective and is bijective if restricted to the subspaces of the corresponding moduli spaces which consist of elements with zero-dimensional stabiliser.
In the classical case $N=2$ we can only have $n=1$ and there are no non-trivial $PU(1)$ connections. Now the components $M^{[F]}$ are the possible contributions to the formula expressing the $PU(N)$ instanton invariant according to the cobordism program indicated above. The generalisation of Witten's conjecture to the $PU(N)$ instanton invariants would follow if only those components $M^{[F]}$ contribute in a non-trivial way for which we have $n=\rk(F)=1$. But the results of \cite{Z_un} indicate that no non-trivial invariants should be expected from $U(n)$ moduli spaces with $n >1$.

The first section sets up our configuration space, introduces the above mentioned quadratic map $\mu_{0,0}$ and derives some important properness property of it. The $PU(N)$ monopole equations are then introduced and it is indicated how to obtain an Uhlenbeck-type compactification of the moduli space. The second sections studies the circle action on the moduli space of $PU(N)$ monopoles, analyses its fixed point set and relates it to $U(n)$ monopoles and $PU(N-n)$ instantons.

\section*{Acknowledgements}
I am very grateful to Andrei Teleman for many mathematical discussions involved with this article. I am also indebted to Peter Kronheimer and Kim Fr\o yshov for discussions on some technical aspects.

\section{The PU(N) - monopole equations}
Here we shall introduce the $PU(N)$-monopole equations associated to the data
of a $Spin^c$-structure $\mathfrak{s}$ and a Hermitian bundle $E \to X$ of rank
$N$ on a Riemannian four-manifold $X$. We shall define the moduli space, and prove a uniform bound on the spinor
component of a solution to the monopole equations.  The standard material in Seiberg-Witten theory ($Spin^c$ structures, $Spin^c$ connections etc.) can be found in one of the textbooks on the topics like \cite{N}, \cite{M} or diverse lecture notes like \cite{T3}.

\subsection{The configuration space}
Let $X$ be a closed oriented Riemannian four-manifold with a
$Spin^c$ structure $\mathfrak{s}$ on it. The $Spin^c$ structure consists of two
Hermitian rank 2 vector bundles $S^\pm_{\mathfrak{s}}$ with identified
determinant line bundles and a Clifford multiplication 
\begin{equation*} 
\gamma : \Lambda^1(T^*X) \to \Hom_\C(S^+_\mathfrak{s},S^-_\mathfrak{s}) \ .
\end{equation*}

Furthermore suppose we are given a
Hermitian vector bundle $E$ with determinant line bundle $w = \det(E)$ on $X$.
We can then form spinor bundles 
\begin{equation*}
W_{\data}^\pm := S^\pm_{\mathfrak{s}} \tensor E .
\end{equation*}
Clifford multiplication extends by tensoring with the identity on $E$. 
This way we obtain a $Spin^c$ - structure `twisted' by the hermitian bundle $E$.
\\

Let $\theta$ be a fixed smooth unitary connection in the determinant line bundle $w$. 
We shall denote by $\mathscr{A}_\theta(E)$ the space of smooth unitary connections on
$E$ which induce the fixed connection $\theta$ in $w$. This is an affine space modelled on
$\Omega^1(X;\mathfrak{su}(E))$. Here $\su(E)$ denotes the bundle of skew-adjoint
trace-free endomorphisms of $E$. Furthermore $\Gamma(X;W^+_{\data})$ denotes the space of
smooth sections of the spinor bundle $W^+_{\data}$.
We define our configuration space to be
\begin{equation*}
\mathscr{C}_{\data,\theta} := \Gamma(X;W^+_{\data}) \times \mathscr{A}_\theta(E) \ .
\end{equation*}
We denote by  $\G$ the group of unitary automorphisms of $E$ with determinant 1; it is the `gauge group' of
our problem. It acts
 in a canoncial way on sections of the spinor bundles, and as $(u,\nabla_A)
\mapsto u \nabla_A u^{-1}$ on the connections, where $u$ is a gauge
transformation and $\nabla_A$ a unitary connection. In particular it lets the induced connection in the
determinant line bundle $w$ fixed. 
The set $\mathscr{B}_{\data}$ is defined to be the configuration space up to gauge, that
is the quotient space $\conf_{\data,\theta} / \G$. 
\\

The reason we consider only smooth objects is purely a matter of simplicity
here. Obviously, as soon as we wish to consider more analytical properties like
transverality, we study suitable Sobolev-completions of these
spaces.
\\


\subsection{Algebraic preliminaries}

We shall now recall the definition \cite{Z_un} of the quadratic map $\mu_{0,\tau}: S^+\tensor E \to \su(S^+)\tensor_\R \su(E)$, defined for a real number $\tau \in [0,1]$. For $\tau = 0$ it will
appear in the $PU(n)$ - monopole equations. This map is a straight-forward generalisation of the corresponding map
appearing in the classical $PU(2)$-situation studied by Teleman and Feehan-Leness and as well as the one in the classical (abelian) Seiberg-Witten equations. 

The twisted spinor bundles $W^\pm_{\data}$ are associated bundles of the fibre product of a
$Spin^c$
principal bundle and a $U(n)$-principal bundle on X, with the standard fibre
$\C^2 \tensor \C^n$.
\\

Let us consider the isomorphism
\begin{equation*}
\begin{split}
  (p,q): \ \,  \gl(\C^n) & \to \sl(\C^n) \oplus \C \, \id \\
   	 a  & \mapsto  \left(a - \frac{1}{n} \ \tr(a) \cdot \id, \frac{1}{n} \
\tr(a) \cdot \id\right) \ .
\end{split}
\end{equation*}
Both components $p$ and $q$ are orthogonal projections onto their images. Note
that $\gl(\C^2) \tensor \gl(\C^n)$ and $\gl(\C^2 \tensor \C^n)$ are canonically
isomorphic. We define the orthogonal projections
\begin{equation*}
\begin{split}
  P: \ \gl(\C^2 \tensor \C^n) & \to \sl(\C^2) \tensor \sl(\C^n) \ , \\
  Q: \ \gl(\C^2 \tensor \C^n) & \to \sl(\C^2) \tensor \C \, \id
\end{split}
\end{equation*}
to be the tensor product $( \ )_0 \tensor p$ respectively $( \ )_0 \tensor q$,
with $( \ )_0$ denoting the trace-free part of the endomorphism of the first
factor $\C^2$. 
 
For elements $\Psi, \Phi \in \C^2 \tensor \C^n$ we define
\[
\mu_{0,\tau}(\Psi,\Phi):=P (\Psi \Phi^*) \ + \tau \, Q ( \Psi \Phi^*) ,
\]
 where $(\Psi \Phi^*) \in \gl(\C^2 \tensor \C^n)$ is defined to be the
endomorphism $\Xi \mapsto
\Psi (\Phi,\Xi)$. 

With this notation $\mu_{0,1} (\Psi,\Phi)$ is simply the orthogonal projection
of the endomorphism $\Psi \Phi^* \in \gl(\C^2 \tensor \C^n)$ onto $\sl(\C^2)
\tensor \gl(\C^n)$. We shall also write $\mu_{0,\tau}(\Psi):=
\mu_{0,\tau}(\Psi,\Psi)$ for the associated quadratic map. In the case $n=1$ the
map $\mu_{0,1}(\Psi)$ is the quadratic map in the spinor usually occuring in the
Seiberg-Witten equations \cite{W} \cite{KM}. The proof of the following proposition can be found in \cite{Z_un}. \\
%

\begin{prop}\label{properness}
Suppose $n > 1$. Then the quadratic map $\mu_{0,\tau}$ is uniformly proper. In
other words, there is a positive
constant $c > 0$ such that 
  \begin{equation}\label{propernessconstant}
    \abs{\mu_{0,\tau}(\Psi)} \geq c \abs{\Psi}^2 \ .
  \end{equation}
As a consequence we have the formula
  \begin{equation}
    \left(\mu_{0,\tau}(\Psi) \Psi, \Psi \right) \geq c^2 \abs{\Psi}^4 \ 
  \end{equation}
whenever $\tau \geq 0$.

\end{prop}

By slightly modifying the proof of the latter proposition in \cite{Z_un} we get also the following:
\begin{prop}\label{no zero divisors}
Suppose $N \geq 2$ or $\tau \neq 0$. 
Then the bilinear map $\mu_{0,\tau}$ is `without zero-divisors' in the following sense: If
$\mu_{0,\tau}(\Psi,\Phi) = 0$, then either $\Psi = 0$ or $\Phi = 0$. 
\end{prop}
\qed 

Because of the equivariance property of the map $\mu_{0,\tau}$ we get in
a straightforward way corresponding maps between bundles, giving rise to
\begin{equation*}
\mu_{0,\tau}: W^\pm_{\data}  \times  W^\pm_{\data} \ 
\to \sl(S^\pm_{\s}) \tensor_{\C} \gl(E) \ ,
\end{equation*}
respectively, for the quadratic map,
\begin{equation*}
\mu_{0,\tau}: W^\pm_{\data} \ 
\to \su(S^\pm_{\s}) \tensor_{\R} \u(E) \ .
\end{equation*}
These maps on the bundle level satisfy the corresponding statement in the above
proposition with the same constant $c$. 
\begin{definition}
  If we wish to make precise to which Hermitian bundle $E$ we refer we shall denote the corresponding bundle as an upper-script $\mu_{0,\tau}^e$.
\end{definition}

\subsection{The PU(N)-monopole equations}

The Clifford map $\gamma$ is, up to a universal constant, an isometry of the
cotangent bundle
onto a real form inside $\Hom_\C(S^+_\mathfrak{s},S^-_\mathfrak{s})$ which can
be specified by the Pauli matrices. We extend $\gamma$ to $\End(S^+_\s
\oplus S^-_\s)$ by $ - \gamma^*$ on the negative Spinor bundle. It then
naturally extends to exteriour powers of $T^*X$, and in particular its
restriction to self-dual two-forms is zero on the negative Spinor bundle, and
induces an isomorphism 
\begin{equation*}
\gamma: \Lambda^2_+(T^*X) \stackrel{\cong}{\to} \su(S^+_\mathfrak{s}) \ .
\end{equation*}

Taking the tensor product with the identity on $E$ induces a Clifford
multiplication $\gamma : \Lambda^1(T^*X) \to \Hom(W^+_{\data},W^-_{\data})$.
Let's fix a background $Spin^c$ connection $B$ on $\mathfrak{s}$. By composing
the connection $\nabla_{B} \tensor \nabla_{\hat{A}}$ with the Clifford
multiplication we get a Dirac operator
\begin{equation*}
D_{\hat{A}} := \gamma \circ \left(\nabla_{B} \tensor \nabla_{\hat{A}}\right) :
\Gamma(X;W^\pm_{\data}) \to \Gamma(X;W^\mp_{\data}) \ .
\end{equation*}
This is a self-adjoint first order elliptic operator. We have
oppressed the $Spin^c$ connection $B$ from the notation because it will not be
a variable in our theory.

We are now able to write down the PU(N)-monopole equations associated to the
data $(\mathfrak{s},E)$ consisting of a $Spin^c$-structure $\mathfrak{s}$ and a
unitary bundle $E$ on $X$. These equations read
as follows:
\begin{equation}\label{PUN-equations}
\begin{split}
 D_{\hat{A}} \Psi & = 0 \\
 \gamma(F_{A}^+) - \mu_{0,0}(\Psi) & = 0 \ . 
\end{split}
\end{equation}

The left hand side of the above equations can be seen as a map $\mathscr{F}$ of
the configuration space $\mathscr{C}= \Gamma(X,W^+_{\data}) \times \Atheta$ to
the space
$\Gamma(X,W^-_{\data})\times \Gamma(X,\su(S^+_{\mathfrak{s}}) \tensor \su(E))$.
As such
it satisfies
the equivariance property
\[
  \mathscr{F}(u.(\Psi,\hat{A})) = (u \times \ad_{u}) (\mathscr{F}(\Psi,\hat{A}))
\ .
\]
Therefore it is sensible to define:
\begin{definition}
The moduli space $M_{\mathfrak{s},E,\theta}$ of $PU(N)$ monopoles is defined
to be the solution set of the equations (\ref{PUN-equations}) modulo the gauge-group $\G$:
\begin{equation*}
M_{\data,\theta} := \{ [\Psi,\hata] \in \bonf_{\data} | \mathscr{F}(\Psi,\hata) = 0 \} \ .
\end{equation*}
\end{definition}
There is an elliptic deformation complex associated to a solution $(\Psi,\hata)$ of the $PU(N)$ monopole equations. Its index equals minus the `expected dimension' of the moduli space. The latter is given by the following formula:
\begin{equation*}
\begin{split}
d(\data) := & -2 \,  \langle p_1(\su(E)),[X]\rangle - \,  (N^2-1) (b_2^+(X) - b_1(X) + 1) \\
			& -\frac{N}{4} \, \text{sign}(X) + \langle \, c_1(E)^2 - 2\, c_2(E) + \, c_1(L)c_1(E) + \, \frac{N}{4} c_1(L)^2 , [X]\rangle \ ,
\end{split}
\end{equation*}
where $p_1(\su(E))$ denotes the first Pontryagin class of the bundle $\su(E)$, $b_2^+(X)$ the maximal dimension of a subspace of $H_2(X;\R)$ on which the intersection form of $X$ is positive definite, and $\text{sign}(X)$ the signature of the intersection form. In this formula the expression in the first line of the right hand side is the expected dimension of the moduli space of $PU(N)$ ASD-connections in $E$, and the second line is the index of the Dirac operator $D_{\hata}$. 

\subsection{Uniform bound, compactification} 
The moduli space $M_{\data, \theta}$ turns out to be non-compact in general, but possesses a canonical compactification very analogue to the Uhlenbeck-compactification of instanton moduli spaces \cite{DK}. The main reason is that there is a uniform $C^0$ bound on the spinor part $\Psi$ of $PU(N)$ monopoles $[\Psi,\hata] \in M_{\data,\theta}$. Knowing this, the compactification is fairly standard \cite{DK}, \cite{T2}, \cite{FL3}, and therefore we will keep our exposition very brief on this point. An outline for the $PU(N)$ case can also be found in \cite{Z1}. 
\\

The $C^0$ bound is derived similarly to classical Seiberg-Witten theory \cite{KM} from the Weitzenb\"ock formula for the Dirac-operator $D_\hata$ by making also use of the above Proposition \ref{properness}, see also \cite{Z_un} for a similar computation.
\begin{prop}\label{uniform bound}
There are constants $c, K \in \R$, with $c > 0$, such that for any
monopole $[\Psi,\hat{A}]\in M_{\data}$ we have
a $C^0$ bound:
\begin{equation}\label{apriori-bound}
	\max \abs{\Psi}^2 \, \leq \, \max \left\{ 0, K/c^2 \right\} \ . 
\end{equation}
Here the constant $K$ depends on the Riemannian metric, the fixed background
$Spin^c$ connection whereas the
constant $c$ is universal.
\end{prop}

Let $\mathfrak{s}$ be a $Spin^c$-structure on $X$ and let $E \to X$ be a
unitary bundle on $X$. We denote by $E_{-k}$ a bundle which has first Chern
class $c_1(E_{-k}) = c_1(E)$ and whose second Chern class satisfies 
\[
\langle c_2(E_{-k}), [X] \rangle = \langle c_2(E), [X] \rangle - k \ .
\]
On a four-manifold such a bundle is unique up to isomorphism.

\begin{definition}
An ideal $PU(N)$ monopole associated to the data $(\data,\theta)$ is given by a pair
$([\Psi,\hat{A}], \bf{x})$, where $[\Psi,\hat{A}] \in
M_{\mathfrak{s},E_{-k},\theta}$ is a $PU(N)$ monopole associated to $(\mathfrak{s},E_{-k})$ monopole, and
$\bf{x}$ is an element of the k-th symmetric power $Sym^k(X)$ of X (that is, an
unordered set of k points in X, ${\bf x} = \{ x_1, \dots, x_k\}$). The
curvature density of $([\Psi,\hat{A}],{\bf x})$ is defined to be the measure
\[
\abs{F_\hata}^2 + 8 \pi^2 \sum_{x_i \in {\bf x}} \delta_{x_i} \ .
\]
The set of ideal monopoles associated to the data $(\data,\theta)$ is
\begin{equation}\label{idealmonopoles}
I M_{\data,\theta} := \coprod_{k \geq 0} M_{\mathfrak{s},E_{-k},\theta}
\times Sym^k(X) \ ,
\end{equation}
\end{definition}
which is endowed with a convenient topology \cite{DK}, \cite{T2}. Rougly speaking, in this topology a sequence in the main stratum $[\Psi_n,\hata_n]$ converges to a point $([\Phi,\hat{B}],{\bf x}) \in M_{\mathfrak{s},E_{-k},\theta}
\times Sym^k(X)$ if the sequence of measures $\abs{F_{\hata_n}}^2 vol$ converges to the measure given by $\abs{F_{\hat{B}}}^2 vol + 8 \pi^2 \sum_{x_i \in {\bf x}} \delta_{x_i} $, and if $\Psi_n$ converges to $\Phi$ in the complement of ${\bf x}$ in $X$. The main result is then:
\begin{theorem}(Compactness-Theorem)
The closure of $M_{\data,\theta}$ inside the space of ideal monopoles
$IM_{\data,\theta}$ is
compact.
\end{theorem}
\section{The circle-action and its fixed-point set, relations to U(n)-monopoles}
This section is the core of our considerations. There is a circle-action on the configuration space
modulo gauge which is induced by multiplying spinors with complex numbers of absolute value one.
The fixed-point set of this circle-action obviously contains elements with zero spinor component, and
the other elements are those which have a connection that splits up into the direct sum of two connections
on proper subbundles on $E$ and which have the spinor component being a section of a corresponding
subbundle. The latter fixed-point set are naturally labelled by isomorphism classes of proper subbundles of
$E$. 
We shall describe a way of `parametrising' these fixed point loci by picking a representative vector
bundle for each isomorphism class. Next we restrict our considerations to the intersection of the
fixed-point set with the moduli space of $PU(N)-$ monopoles: fixed-points with vanishing spinor are then
simply anti-selfdual $PU(N)$- connections in $E$, whereas the other fixed-point sets are fibrations of
moduli spaces of $PU(n)$-connections in a summand $F$ of $E$ of rank $n$ over moduli spaces of
$U(N-n)$-monopoles in the complement of $F$ in $E$.

\subsection{Reductions and stabilisers of connections under the gauge group}
Here we study the stabilisers of connections $\hata \in \Atheta$ under
the action of the gauge group $\G$. 


An element $u$ of the gauge group $\G$ acts on the set of connections $\hat{A}$
in $\Atheta$ by the formula 
\[
  u(\hat{A}) = \hat{A} - (d_{\hat{A}} u) u^{-1} \ ,
\]
where we consider $u$ as section of the vector bundle $\gl(P_E)$. The
stabiliser of the connection $\hat{A}$ inside the gauge
group $G$ is the group of automorphisms which preserve $\hat{A}$:
\begin{equation*}
\begin{split}
  \Gamma(\hat{A}) = & \left\{ u \in \G \ | \ u(\hat{A}) = \hat{A} \right\} \\
   = & \left\{ u \in \G \ | \ d_{\hat{A}} u = 0 \right\}
\end{split}
\end{equation*}
This group $\Gamma(\hat{A})$ is a finite-dimensional compact Lie-group if the
base manifold $X$ is compact. Indeed, we can see by parallel transport that an
element $u \in \Gamma(\hat{A})$ is determined by restrictions to one fibre in
each component of $X$. In particular, if $X$ is connected, then
$\Gamma(\hat{A})$ can be seen as a closed Lie subgroup of $\Aut(E_x)$ for any
point $x \in
X$. Further, by choosing an element in the fibre of $x$ inside the principal
fibre bundle $P_E$, we get a non-canonical identification of $\Aut(E_x)$ with
the structure
group $U(N)$. Thus $\Gamma(\hat{A})$ can be seen as a closed Lie supgroup of
$U(N)$.

Note that the stabiliser always contains the centre $Z(G)$ of the
structure group $G$,
$\Gamma(\hat{A}) \supseteq Z(G)$. In our case, the centre injects as
\[
 Z(SU(N)) = \left\{ \lambda \ \id_E \ | \ \lambda^N = 1  \right\}
\]
\begin{definition}
The connection $\hat{A} \in \Atheta$ is called {\em reducible} if the
stabiliser $\Gamma(\hat{A})$ is different from the centre $Z(SU(N))$. 
\end{definition}
Suppose now that $\hat{A}$ is reducible and that $d_{\hat{A}} u = 0$. We recall
a standard result for normal endomorphisms of a Hermitian vector-bundle. An
endomorphism $u$ is called normal if $u u^* = u^* u$, where $u^*$ is the adjoint
endomorphism with respect to the Hermitian structure. The following lemma is easy to prove:
\begin{lemma}\label{decomposition_of_parallel_endomorphism}
  Suppose the normal endomorphism $u \in \End(E)$ is $\hat{A}$-parallel, $d_\hata u = 0$. Then
its spectrum is constant and there is a $\hat{A}$-parallel decomposition of $E$
into subbundles
  \begin{equation*}
    E = \bigoplus_{\lambda \in \text{Spec}(u)} E_{\lambda} \ .
  \end{equation*}
Each summand $E_\lambda$ is $u$-invariant, and we
have $u|_{E_\lambda} = \lambda \ \id_{E_\lambda}$. In other words, the
$E_\lambda$ are eigen-bundles of the endomorphism $u$.
\end{lemma}
 \qed
As a corollary one obtains that a connection $\hata$ is reducible if and only if there is a proper
subbundle $F$ of $E$ which is $\hata$ - parallel. For, the latter condition clearly implies that $\hata$ is
reducible in our definition. On the other hand, if the stabiliser $\Gamma(\hata)$ is strictly bigger than
the centre $Z(SU(N))$, then by the preceding lemma there must be a $\hata$-parallel automorphism $u \in
\G$ which admits an eigenvalue which is not an $N^{th}$ root of one, and therefore there must be such a
$\hata$-parallel subbundle.

\subsection{Stabiliser of a configuration under the gauge group}
The topology of the configuration space up to gauge $\bonf_{\data} =
\conf_{\data} / \G$ does not really have nice properties. It can be quite
singular due to the fact that the gauge-group may have non-trivial stabilisers.
However, if we restrict our attention to the subset of the configurations which
have trivial stabilisers, then the quotient under the gauge-group $\G$ has the
nice property of being a Banach-manifold \footnote{actually, we rather get a
`Fr\'echet'-manifold, since we are working with $C^\infty$ configurations - but
as soon as one completes the configuration space with respect to some Sobolev
norms we really get a Banach-manifold}. If we require the stabilisers to be
finite groups, then we still get a quotient in which the singularities are
relatively mild. However, the configurations which have positive-dimensional
stabilisers are more delicate.\\

We define the stabiliser $\Gamma(\Psi,\hat{A})$ of a configuration
$(\Psi,\hat{A}) \in \mathscr{C}_{\data}$ to be the set 
\begin{equation*}
\begin{split}
  \Gamma((\Psi,\hat{A})) = & \left\{ u \in \G \ | \ u(\Psi,\hat{A}) =
(\Psi,\hat{A})
\right\} \\
   = & \left\{ u \in \G \ | \ u(\Psi) = \Psi, \quad d_{\hat{A}} u = 0 \right\} \
.
\end{split}
\end{equation*}

\begin{definition}
  The subset $\mathscr{C}_{\data}^* \subseteq \mathscr{C}_{\data}$
(respectively $\mathscr{C}_{\data}^{**}$) is
defined to be the set of configurations $(\Psi, \hat{A}) \in
\mathscr{C}_{\data})$ which has zero-dimensional stabiliser (respectively
trivial stabiliser). The subset
$\mathscr{B}_{\data}^* \subseteq \mathscr{B}_{\data}$ is defined to be the
subset $\mathscr{C}_{\data}^* / \G \subseteq \mathscr{B}_{\data}$. The subset
$\mathscr{B}_{\data}^{**}$ is defined correspondingly.
\end{definition}

It is easy to see that the stabiliser of a configuration $(\Psi,\hat{A})$, with
$\Psi$ non-vanishing and the connection $\hat{A}$ irreducible, is trivial. 
Furthermore, the configurations $(\Psi,\hat{A})$ with $\Psi \equiv 0$ and
irreducible connection part $\hat{A}$ have stabilisers which are the finite
group $\Z_N$. These claims, as well as the following Proposition, follow easily
from Lemma (\ref{decomposition_of_parallel_endomorphism}) in the section of
reducible connections above.

\begin{prop}\label{criterion_finite_stabiliser}
  Suppose we have given a configuration $(\Psi,\hat{A})$ with non-vanishing 
spinor $\Psi$, and the connection of the form $\hat{A} = \hat{A}_1 \oplus
\hat{A}_2$ according to a $\hat{A}$-parallel decomposition $E=F \oplus F^\perp$,
with both $\hat{A}_1$ and $\hat{A}_2$ irreducible. Then its stabiliser is a
finite (and thus zero-dimensional) group.
\end{prop}

\subsection{The circle action}
We are given an $S^1$-action on the configuration space $\mathscr{C}_{\data}$
given by
the simple formula
\begin{equation*}
\begin{split}
  S^1 \times \mathscr{C}_{\data} \to & \mathscr{C}_{\data} \\
  \left(z,(\Psi,\hat{A})\right) \mapsto & \left(z \Psi, \hat{A}\right)
\end{split}
\end{equation*}
Now as this action commutes with the action of the gauge group $\G$, we see that
we get a well-defined action on the quotient,
\[
S^1 \times \bonf_{\data,\theta} \to \bonf_{\data,\theta} 
\]
The action is not effective. In fact if $z^N=1$ there is always a
gauge-transformation $u$ with $u(\Psi,\hat{A}) = (z \Psi, \hat{A})$, that is
$[\Psi,\hat{A}] = [z \Psi, \hat{A}]$. This is because the centre $Z(SU(N)) \cong \Z/N$ always injects into
the stabiliser of the connection $\hata \in \Atheta$. Therefore we define
\[
   r\left(z, \left[\Psi, \hat{A}\right]\right) := \left[ z^{1/N} \Psi, \hat{A}
\right] \ .
\]
In this formula $z^{1/N}$ is any $N^{th}$ root of $1$, the equvialence class $[ z^{1/N} \Psi, \hat{A}]$ does
not depend on the particular choice.

\begin{remark}
Suppose we had chosen as gauge-group $\G$ the group of unitary
bundle-automorphisms of $E$ which fix the connection $\theta$ only, that is,
the larger group of unitary automorphisms with constant determinant. Then 
the same action on $\mathscr{C}_{\data} = \Gamma(X,W^+) \times \Atheta$ would
 have introduced the trivial $S^1$-action on $\mathscr{B}_{\data}$. This
justifies our choice of the gauge group $\G$ as $\Gamma(SU(E))$.
\end{remark}

\begin{lemma}
\label{fixpoints}
Suppose there exists some $z_0 \in S^1$ with $z_0\neq 1$ such that
$[z_0^{1/N}\Psi,\hat{A}] = [\Psi, \hat{A}]$. Then for any $z \in S^1$ we have 
  \[
    [z\Psi,\hat{A}] = [\Psi, \hat{A}] \ .
  \]
Thus an element $[\Psi,\hata]$ is a fixed-point of the circle-action if and only if it is left fixed by some
non-trivial element $z_0 \in S^1$ under the action $r$.
\end{lemma}
{\em Proof:}
This is certainly true for vanishing spinor component. Therefore suppose $\Psi \neq 0$. By hypothesis there
is some $z_0\neq 1$, a gauge-transformation $u \in \G$ and an $N^{th}$ root $z_0^{1/N}$ of $z_0$ such that 
\begin{equation*}
   u (\Psi)  =  z_0^{1/N} \Psi  \ ,
\end{equation*}
and furthermore we have $d_\hata u = 0$. Therefore $z_0^{1/N}$ is an eigenvalue of the endomorphism $\id \tensor u$ on $S^+ \tensor E$, but this implies that $z_0^{1/N}$ is an eigenvalue of $u$ which is not an $N^{th}$ root of the unity. By the above lemma \ref{decomposition_of_parallel_endomorphism} we see that $\hata$ must be reducible, $E$ splits into a direct sum of $\hata$-parallel summands $E = \oplus E_i$ and $u$ decomposes into $u = \sum \lambda_i \id_{E_i}$. One of the eigenvalues $\lambda_i$ must be equal to $z_0^{1/N}$ and $\Psi$ is a section of the corresponding bundle $S^+ \tensor E_i$. With this given splitting and for any $z \in S^1$ one then can certainly find gauge transformations $u_{z^{1/N}}$ for any $N^{th}$ root $z^{1/N}$ of $z$ such that $d_\hata u_{z^{1/N}} = 0$ and $u_{z^{1/N}} \Psi = z^{1/N} \Psi$.  \qed

Using the above lemma \ref{decomposition_of_parallel_endomorphism} again we deduce the following simple criterium for fixed-points under the action $r$ above:

\begin{prop}\label{fixed points}
A configuration up to gauge $[\Psi,\hat{A}] \in \mathscr{B}_{\data,\theta}$ is
contained in the fixed point set of the action $r$ if and only if for some (or
equivalently, for any) representative $(\Psi,\hat{A})$ we have one of the
following (possibly both):
\begin{enumerate}
\item There is a non-trivial $\hat{A}$-parallel orthogonal decomoposition $E=
\oplus 
E_i$ and the spinor is a section of one of $S^+_\mathfrak{s} \tensor E_i$
\item The spinor vanishes $\Psi \equiv 0$ \ .
\end{enumerate}
\end{prop}
Further down we will see that if we impose in addition the monopole equations the spinor component of a fixed point $[\Psi,\hata]$ will automatically lie in a proper summand $S^+ \tensor E_i$ as soon as the connection $\hata$ is reducible.

\subsection{The $S^1$-fixed point set inside the
configuration space modulo gauge, and parametrisations}

First we will describe the fixed-point set of the $S^1$-action inside
$\bonf_{\data}$. In the above Proposition \ref{fixed points} we saw that these are related to proper subbundles of $E$. However, two subbundles which are mapped into each other by gauge transformations, i.e. automorphisms of E, should be considered equivalent. This equivalence of subbundles might be called `ambiently isomorphic', but it is easy to see that two subbundles of $E$ are ambiently isomorphic if and only if they are isomorphic as abstract bundles. This even holds for prescribed determinant. For the further work, especially for describing the
intersection of the fixed point set with the moduli space, it turns out useful
to fix representatives $F \in [F]$ for each such isomorphism class. This 
yields to a `parametrisation' of each component of the fixed point set which is
determined by the isomorphism class $[F]$. \\

\begin{definition}\label{def_klammerf}
Let $F$ be a proper summand of the unitary bundle $E$. We define the set
$\mathscr{B}_{\data}^{[F]}$ to be the set of all elements $[\Psi,\hat{A}]
\in \mathscr{B}_{\data}$ such that for some
representative $(\Psi,\hat{A})$ there exists a $\hat{A}$-parallel
decomposition $E=F \oplus F^{\perp}$ with $\Psi \in
\Gamma(X,S^+_{\mathfrak{s}} \tensor F)$.
\end{definition}
It is easy to see that this set $\mathscr{B}_{\data}^{[F]}$ is contained in the
$S^1$-fixed point set $\mathscr{B}_{\data}^{S^1}$. Another subset of the fixed
point set is given by the following definition:

\begin{definition}
We define the set $\mathscr{B}_{\data}^{\equiv 0}$ to be the set of set of all
elements $[\Psi,\hat{A}]
\in \mathscr{B}_{\data}$ such that for some (or equivalently, for any)
representative $(\Psi,\hat{A})$ the spinor vanishes identically, $\Psi \equiv
0$. This set is contained in the $S^1$-fixed point set.
\end{definition}


The above Proposition \ref{fixed points} gives then the following description of
the fixed point set of the $S^1$ - action $r$:
\begin{prop}
The $S^1$-fixed point set $\mathscr{B}_{\data}^{S^1}$ is given as the union
\[
 \left( \bigcup_{[F] \subseteq E} \mathscr{B}_{\data}^{[F]} \right) \ \bigcup \ 
\mathscr{B}_{\data}^{\equiv 0} \ .
\]
Here the first union is taken over all isomorphism classes of
proper subbundles of $E$.
\end{prop}
We should point out as well that the different components $\bonf_{\data}^{[F]}$
may a priori intersect each other or the fixed-point component of vanishing
spinor $\bonf_{\data}^{\equiv 0}$. We would also like to remark that for rank strictly higher than $2$ we may always have
infinitely many such isomorphism classes of proper subbundles $[F]$ of $E$, even for definite intersection form. 

%

In order to have a convenient description of the set $\bonf_{\data,\theta}^{[F]}$ it seems natural to fix an actual proper subbundle $F$ for each isomorphism class $[F]$. 
Hence the following definition:

\begin{definition}
We define the configuration space relative to the splitting $E=F \oplus
F^\perp$ as the following set:
\begin{equation*}
\begin{split}
\mathscr{C}_{\data, \theta}^{F \oplus F^{\perp}} := & \left\{  (\Psi,\hat{A}_1, 
\hat{A}_2) \in \Gamma(S^+_\mathfrak{s}\tensor F) \times \mathscr{A}(F)
\times \mathscr{A}(F^\perp) \right| \\ &
  \left. \det(\hat{A}_1) \tensor \det(\hat{A}_2)
= \theta \right\}
\end{split}
\end{equation*}
Correspondingly, the group of unitary automorphisms with determinant 1
respecting the splitting $E= F \oplus F^\perp$ is defined to be
\begin{equation*}
\begin{split}
\mathscr{G}^0_{F \oplus F^\perp}:= \left\{ \left. (u_1,u_2) \in \Gamma(U(F))
\times \Gamma(U(F^\perp)) \right| \det(u_1) \cdot \det(u_2) = 1 
 \right\}
\end{split}
\end{equation*}
As usually, we denote the quotient by: 
\begin{equation*}
 \mathscr{B}_{\data,\theta}^{F \oplus F^\perp}:= \mathscr{C}_{\data,\theta}^{F \oplus
F^{\perp}} / \mathscr{G}^0_{F \oplus F^\perp} \ .
\end{equation*}
\end{definition}
It is then easy to see that this yields a well-defined map 
\begin{equation*}
\begin{split}
i_F: &  \mathscr{B}_{\data}^{F \oplus F^\perp} \to \mathscr{B}_{\data}^{[F]}
\\
 	& [ \Psi, \hat{A}_1, \hat{A}_2 ] \mapsto [\Psi, \hat{A}_1 \oplus
\hat{A}_2 ] \ .
\end{split}
\end{equation*}
This map is easily seen to be always surjective. However, it fails to be
injective in general. Nonetheless, on a dense subset of $\mathscr{B}_{\data}^{F
\oplus F^\perp}$ it is, as we shall show next. We will think of the map $i_F$
as a `parametrisation' of the fixed-point set component $\bonf_{\data}^{[F]}$.

As before, let us denote by
$\mathscr{B}^{* \ F \oplus F^\perp}_{\data}$ and by 
$\mathscr{B}^{* \ [F]}_{\data}$ the configurations which have
finite-dimensional
stabiliser in their groups $\G$ respectively $\mathscr{G}^0_{F \oplus F^\perp}$.
Further, we denote by $\mathscr{B}^{*,irr \ F \oplus F^\perp}_{\data}$ the
subset of $\mathscr{B}^{* \ F \oplus F^\perp}_{\data}$ consisting of
elements $[\Psi,\hat{A}_1,\hat{A}_2]$ with non-vanishing spinor, $\Psi \neq 0$,
and both connections $\hat{A}_1$ and $\hat{A}_2$ irreduible. By the
way, $(\Psi,\hat{A}_1,\hat{A}_2)$ has zero-dimensional stabiliser in $\G_{F
\oplus F^\perp}$ if and only if $(\Psi,\hat{A}_1 \oplus \hat{A}_2)$ does so in
$\G$. Now we can state the following:

\begin{prop}\label{Prop_of_iF}
Restriction of the map $i_F$ yields an injective map
\begin{equation*}
i_F: \  \mathscr{B}^{*,irr \ F \oplus F^\perp}_{\data}
\to 
\ \mathscr{B}^{* \ [F]}_{\data}
\end{equation*}
from the subset $\mathscr{B}^{*,irr \ F \oplus F^\perp}_{\data}$ of
configurations, up to gauge, with zero-dimensional stabilisers inside
$\G_{F \oplus F^\perp}$ and irreducible connections, into the fixed-point set
component $\mathscr{B}^{* \ [F]}_{\data}$ of configurations, up to gauge, with
zero-dimensional stabilisers inside $\G$.
\end{prop}
{\em Proof:}
For simplicity we note $i$ instead of $i_F$. Suppose we have elements
$[\Psi,\hat{A}_1, \hat{A}_2]$, $ [\Phi,\hat{B}_1,\hat{B}_2] \in
\mathscr{B}_{\data}^{*,irr \, F \oplus F^\perp} $ such that $i([\Psi,\hat{A}_1,
\hat{A}_2] = i [\Phi,\hat{B}_1,\hat{B}_2]$. This is equivalent to saying that
there is a gauge transformation $u \in \G$ such that 
\begin{equation*}
\begin{split}
u (\Psi) & = \Phi \\
u (\hat{A}_1 \oplus \hat{A}_2 ) & = \hat{B}_1 \oplus \hat{B}_2 \ .
\end{split}
\end{equation*}
The second equation implies that $u$ is an $(\hat{A}_1 \oplus \hat{A}_2)
\tensor ( \hat{B}_1 \oplus \hat{B}_2)^*$-parallel endomorphism of $E$. Let us
write $u$ in the form 
\[
\begin{pmatrix} u_{11} & u_{12} \\ u_{21} & u_{22} \end{pmatrix}
\]
according to the splitting $E= F \oplus F^{\perp}$. Injectivity will follow if
we have $u_{12}= 0 $ and $u_{21} = 0$. It is enough to show just $u_{21}=0$,
as the other equation will follow from the fact that $u$ is unitary. We find
that the morphism $u_{11}$ is $\hat{B}_1 \tensor \hat{A}_1^*$-parallel,
the morphism $u_{12}$ is $\hat{B}_1 \tensor \hat{A}_2^*$-parallel, 
the morphism $u_{21}$ is $\hat{B}_2 \tensor \hat{A}_1^*$-parallel, and 
the morphism $u_{22}$ is $\hat{B}_2 \tensor \hat{A}_2^*$-parallel . 

Now all the connections $\hat{A}_{i}$, $\hat{B}_{i}$, $i=1,2$ are {unitary}
connections. Therefore the fact that, for instance, $u_{12}$ is $\hat{B}_1
\tensor \hat{A}_2^*$-parallel implies that the adjoint $u_{12}^*$ is $\hat{A}_2
\tensor \hat{B}_1^*$-parallel. As a consequence,
the endomorphism $u_{21}^* u_{21}$ of $F$ is $\hat{A}_1 \tensor
\hat{A}_1^*$-parallel, and 
the endomorphism $u_{12}^* u_{12}$ of $F^\perp$ is $\hat{A}_2 \tensor
\hat{A}_2^*$-parallel. 
By the hypothesis $\hat{A}_1$ and $\hat{A}_2$ are irreducible, so that the
above Lemma \ref{decomposition_of_parallel_endomorphism} implies that there are
constants $\xi, \zeta \in \C$ with 
\begin{equation*}
\begin{split}
u_{21}^* u_{21} & = \xi \ \id_{F} \\
u_{12}^* u_{12} & = \zeta \ \id_{F^\perp} \ .
\end{split}
\end{equation*}
We have to show now that under our hypothesis $\xi = 0 $ or $\zeta = 0$,
implying then that $u_{21}= 0$ respectively $u_{12}=0$. But if we had $\xi \neq
 0$, then $u_{21}$ is injective at each point $x \in X$. By the hypothesis we
get that 
$\Psi \neq 0$, and therefore we would have a non-trivial section $u_{21}
(\Psi) \in
S^+_\mathfrak{s} \tensor F^\perp$. However, we have $u(\Psi) = \Phi$, where
$\Phi$ is a section of $S^+_\mathfrak{s} \tensor F$, so that this would yield a
contradiction. Therefore $\xi = 0$ and as a consequence $u_{21}= 0$ and $u_{12}
= 0$. 
\qed 
\begin{remark}
From the above proof we see that $\bonf_{\data}^{*,irr \ F \oplus F^\perp}$ is
not
the maximal possible subset on which $i_{F}$ is injective. For instance, it
would have been enough that $\Psi \neq 0$ and
$\hat{A}_1$ irreducible, but then we do not necessarily have value inside the
configurations with zero-dimensional stabilisers. 
\end{remark}

Next we shall discuss a canonical fibering of the configuration space up to
gauge respecting the proper decomposition $E = F \oplus F^\perp$ that we have
introduced above. Let us denote now by $\G_{F}$ the group of special unitary automorphisms of the
unitary bundle $F$ on $X$, that is $\G_F = \Gamma(X,SU(F))$. So, with this
 notation, $\G_E$ is the
gauge group we have until now denoted by the letter $\G$. On the other hand, we
shall denote by $\mathscr{G}_{F}$ the group of unitary automorphisms of $F$,
that is, $\mathscr{G}_{F} = \Gamma(X,U(F))$. 
\begin{lemma}
We have an exact sequence of groups given by
\begin{equation*}
1 \to \G_{F^\perp} \stackrel{i}{\to} \G_{F \oplus F^\perp} \stackrel{j}{\to}
\mathscr{G}_{F} \to 1 \ . 
\end{equation*}
Here the morphisms are given by $i(u_2) = (\id_F, u_2)$ and $j((u_1,u_2)) :=
u_1$.
\end{lemma}
{\em Proof:} The only non-trivial point is the surjectivity of the morphism
$j$. Indeed, for a given gauge transformation $u_2 \in \mathscr{G}_{F}$ we have
to find some automorphism $u_1 \in \mathscr{G}_{F^\perp}$ such that $\det(u_1)
\cdot \det(u_2) = 1$. So we have to find an automorphism of $F^\perp$ with
prescribed determinant $\det(u_2)^{-1}$. That this is indeed possible follows
from obstruction theory \cite{S} \cite{MS}. \qed

We shall introduce some new notation now. Given a Hermitian vector bundle $F$ on $X$ we shall denote by $P_F$ its associated frame bundle, a principal bundle of structure group $U(n)$, where $n$ is the rank of $F$. Let us denote
by $\mathscr{A}_{PU}(F)$ the affine space of connections in the associated 
$PU(n)$-bundle $P_{F}\times_\pi PU(n)$, where $\pi$ is the natural projection $U(n)\to PU(n)$.
Note that in the case that $n= \text{rank}(F) = 1$ the bundle
$P_{F}\times_\pi PU(n)$ is the trivial principal bundle with structure
group the trivial group, and both $\mathscr{A}_{PU}(F)$ and
$\mathscr{A}_{PU}(F) / \G_{F}$ consist of a single point.
\begin{definition}
We introduce the configuration space
$\conf_{\dataf}^{U}$ to be the set of all
configurations
$(\Psi,\hat{B}) \in \Gamma(S^+_\mathfrak{s} \tensor F) \times \mathscr{A}(F)$.
The subscript $U$ points
out that here we take arbitrary unitary connections in the bundle $F$ (and
not such with fixed determinant). Correspondingly, $
\mathscr{B}^U_{\mathfrak{s},F}$ denotes the configuration space
$\conf_{\dataf}^{U}$ quotiented by the gauge-group $\mathscr{G}_F$ of full
unitary automorphisms of the unitary bundle $F$.
We shall also denote by $\mathscr{B}_{F}^{PU}$ the set of all $PU(n)$-connections $A \in
\mathscr{A}_F^{PU}$ in the
unitary bundle $F$ up to the gauge group $\G_{F}$ of special unitary
automorphisms of the bundle $F$.
\end{definition}

\begin{prop}\label{map_h}
Suppose the 4-manifold $X$ is simply connected. Then we have a bijection
\begin{equation*}
\begin{split}
h: & \mathscr{B}_{\data}^{F \oplus F^\perp}  \stackrel{\cong}{\to}
\mathscr{B}^U_{\mathfrak{s},F} \times \mathscr{B}^{PU}_{F^\perp}
 \\
& [ \Psi,\hat{A}_1, \hat{A}_2 ]  \mapsto ([\Psi,\hat{A}_1],[A_2]) \ .
\end{split}
\end{equation*}
\end{prop}
{\em Proof:} 

First of all it is easily checked that the above map $h$ is well-defined. Let $n$ denote the rank of $F$, $n=\rk(F)$. The rank of $F^\perp$ is then $N-n$. 
At
this point we should recall the action of the gauge-groups
$\mathscr{G}_{F^\perp}$ and $\G_{F^\perp}$ on $PU(N-n)$-connections in the
associated bundle $P_{F^{\perp}} \times_\pi PU(N-n)$. It suffices to discuss
the first case. So let $u_2$ be a unitary automorphism of $F^\perp$, or,
equivalently, an automorphism of the $U(N-n)$-bundle $P_{F^\perp}$. It induces
an automorphism $\bar{u}_2$ of the associated bundle $P_{F^{\perp}} \times_\pi
PU(N-n)$ by the formula $\bar{u}_2 ([p,a]) := [u_2(p),a]$, where here $p \in
P_{F^\perp}$ and $a \in PU(N-n)$, and $[p,a]$ denoting the associated element in
the associated bundle. Let us denote by $\phi_{u_2}$ the $\Ad$-equivariant map
$P_{F^\perp} \to U(N-n)$ associated to $u_2$, and $\phi_{\bar{u}_2}$ the
$\Ad$-equivariant map $P_{F^\perp} \times_\pi PU(N-n)$ associated to $u_2$. The
two are related by the formula
\[
\pi \circ \phi_{u_2} = \phi_{\bar{u}_2} \circ \bar{\pi} \ ,
\]
where $\bar{\pi}$ is the natural bundle morphism $P_{F^\perp} \to P_{F^\perp}
\times_\pi PU(N-n)$ given by $\bar{\pi}(p):=[p,1]$. The connection $\bar{u}_2
(A_2)$ is given by the formula:
\[
\bar{u}_2 (A_2) = A_2 - \phi_{\bar{u}_2}^{-1} d_{A_2} \phi_{\bar{u}_2} \ .
\]
Here $d_{A_2}$ denotes the covariant derivative associated to the connection
$A_2$. On the other hand, if $A_2$ is the $PU(N-n)$-connection induced
 from the unitary connection $\hat{A}_2 \in \mathscr{A}(F^\perp)$,
then $\bar{u}_2 (A_2)$ is the $PU(N-n)$-connection induced from the unitary
connection $u_2 (\hat{A}_2)$. Therefore we shall also denote
$\bar{u}_2(A_2)$ simply by $u_2(A_2)$, bearing in mind, however, that $u_2
(A_2)$ does not depend on $\det(u_2)$. This means that whenever 
$u_2$ and $u_2'$ differ by a $U(1)$-valued function $\zeta$, then the
automorphism $u_2'= \zeta u_2$ has the same effect on $A_2$, $u_2'(A_2) =
u_2(A_2)$. With these preliminaries by hand we can now succeed with the proof. 

We would like to construct an inverse of the map $h$ introduced in the
Proposition we want to prove. First we should have an idea what it should look
like on the configuration space level (before equivalence up to gauge). So
the entity $((\Psi,\hat{A}_1),A_2)$ should be mapped to a convenient element in
$\mathscr{C}_{\data}^{F \oplus F^\perp}$. Now let us denote by
$\hat{B}_2(\theta,A_2,\hat{A}_1)$ the unique unitary connection $\hat{B}_2 \in
\mathscr{A}(F^\perp)$ such that 
\begin{enumerate}
\item \ $ \det(\hat{A}_1) \tensor \det(\hat{B}_2) = \theta \quad \Leftrightarrow
\quad \det(\hat{B}_2) = \theta \tensor \det(\hat{A}_1^*)$  \ , and
\item \ The $PU(N-n)$-connection induced from $\hat{B}_2$ is $A_2 \in
\mathscr{A}_{PU}(F^\perp)$ .
\end{enumerate}
So on the configuration space level we have the map
\begin{equation*}
\begin{split}
 k: &  \mathscr{C}_{\mathfrak{s},F}^U \times \mathscr{A}_{PU}(F^\perp) 
  \to \mathscr{C}_{\data}^{F \oplus F^\perp} \\
 & \left((\Psi,\hat{A}_1),A_2\right) \mapsto
\left(\Psi,\hat{A}_1,\hat{B}_2(\theta,A_2,\hat{A}_1) \right) \ . 
\end{split}
\end{equation*}
We would like to show that it descends to the quotients by the respective
gauge-groups. So let $u_2 \in \G_{F^\perp} = \Gamma(X,SU(F^\perp))$ and $u_1
\in \mathscr{G}_{F} = \Gamma(X,U(F))$. Then the entity 
$(u_1(\Psi,\hat{A}_1),u_2(A_2))$ maps to 
\[
\left(u_1 \Psi, u_1(\hat{A}_1),
\hat{B}_2(\theta,u_2(A_2),u_1(\hat{A}_1))\right). 
\]
We have to show that there is a gauge transformation $(u_1',u_2') \in \G_{F
\oplus F^\perp}$, that is, unitary automorphisms $u_1',u_2'$ with $\det(u_1')
\cdot \det(u_2') = 1$, such that 
\begin{equation*}
(u_1',u_2'). \left(\Psi,\hat{A}_1, \hat{B}_2(\theta,A_2,\hat{A}_1) \right) = 
\left(u_1 \Psi, u_1(\hat{A}_1),
\hat{B}_2(\theta,u_2(A_2),u_1(\hat{A}_1))\right) \ .
\end{equation*}
Clearly we must have $u_1'= u_1$, at least if the spinor is non-vanishing, but
otherwise it appears to be the natural choice. Now by assumption $X$ is simply
connected, so that any $U(1)$-valued function on $X$ admits an $(N-n)$-th root,
albeit not a canonical one. 
%
%
%
%
Let us denote by
$\det(u_1)^{-1/(N-n)}$ {\em some} choice of an $(N-n)$-th root of $\det(u_1)^{-1}$.
With the preliminaries made, and by explicit computations, it
is then easy to see that $u_2' := \det(u_1)^{-1/(N-n)} u_2$ will satisfy the requirement. This
shows that the map $k$ induces a map $\bar{k}$ when passing to the quotients:
\begin{equation*}
\begin{split}
\bar{k}:
 \mathscr{B}^U_{\mathfrak{s},F} \times \mathscr{B}_{PU}(F^\perp) \to
\mathscr{B}_{\data}^{F \oplus F^\perp} \ .
\end{split}
\end{equation*}
It is now easily checked that $\bar{k}$ and $h$ are inverses of each other. 
\qed
\begin{remark}
  Without the assumption that $X$ is simply-connected we can still show that we
get a fibration $\mathscr{B}_{\data}^{F \oplus F^\perp} \to
\mathscr{B}_{\mathfrak{s},F}^U $ with standard fibre
$\mathscr{A}_{PU}(F^\perp)/\mathscr{G}^0_{F^\perp}$. The non-triviality of this
fibration should be encoded in $H_1(X,\Z)$. From now on, however, we shall
suppose that our 4-manifold $X$ is simply connected.
\end{remark}

\subsection{The circle-action on the moduli space of $PU(N)$-monopoles}

Until now our consideration of the $S^1$-action and its fixed point set was
inside the configuration space up to gauge, $\mathscr{B}_{\data,\theta}$. 
It is an easy observation that the moduli space of $PU(N)$-monopoles 
$M_{\data,\theta}\subseteq \mathscr{B}_{\data,\theta}$ is invariant under the
$r$-action, $r(S^1,M_{\data}) \subseteq M_{\data}$. 
All we have found out about the circle-action on $\bonf_{\data,\theta}$ applies to the restriction of this action to the moduli space as well. However, there are more things we can say about the fixed-point set of the circle-action for this restriction. In particular, these fixed-point sets are naturally related to other moduli spaces. Obviously the intersection $\bonf_{\data,\theta}^{\equiv 0} \cap M_{\data,\theta}$ consists of anti-self-dual connections in $E$, and the intersection $\bonf_{\data,\theta}^{[F]} \cap M_{\data,\theta}$ is parametrised by the product of the moduli space of $U(n)$-monopoles in $F$, with $n=\rk(F)$, and the moduli space of anti-self-dual $PU(N-n)$-connections in $F^\perp$.  

Proposition \ref{fixed points} above described the fixed-points of the circle-action on $\bonf_{\data,\theta}$, the configuration space modulo gauge. In particular the element $[\Psi,\hata]$ lies in $\bonf_{\data,\theta}^{[F]}$ if and only if for a representative $(\Psi,\hata)$ we have a $\hata$-parallel decomposition $E= F \oplus F^\perp$, and the spinor part $\Psi$ is a section of $S^+_\mathfrak{s} \tensor F$. This second condition becomes automatically satisfied if $(\Psi,\hata)$ solve the $PU(N)$ monopole equations: 

\begin{prop}\label{one_component_vanishes}
Suppose the configuration $(\Psi,\hat{A})$ satisfies the
$PU(N)$-Seiberg-Witten-equations (\ref{PUN-equations}) associated to the data
$(\data)$. \\ Sup\-pose further that the
connection $\hat{A}$ is reducible, and that $E = \oplus E_{i} $  is a
$\hat{A}$-parallel orthogonal decomposition into proper subbundles, and that the
base manifold $X$ is
connected. Then the spinor must be a section of one of the
bundles $S^+_\mathfrak{s} \tensor
E_i$. 
\end{prop}
{\em Proof:} Suppose the connection $\hat{A}$ splits into two connections
$\hat{A_1}
\oplus \hat{A_2}$ with respect to $E = E_1 \oplus E_2$. As an endomorphism of
$E$ the curvature $F_{\hat{A}}$ splits as 
\begin{equation*}
  F_{\hat{A}} = \begin{pmatrix} F_{\hat{A_1}} & 0 \\ 0 & F_{\hat{A_2}} \ 
\end{pmatrix} \ .
\end{equation*}
In other words, it is a section of $\Lambda^2(T^*X) \tensor \left(\u(E_1) \oplus
\u(E_2)\right)$. The trace-free part $F_A$ is then a section of the bundle
$\Lambda^2(T^*X) \tensor ((\u(E_1) \oplus \u(E_2))\cap \su(E))$. Therefore the
curvature-equation of the $PU(N)$-monopole-equations implies that 
\begin{equation}\label{mu-decomposition}
\mu_{0,0}(\Psi) \in \Gamma(X,\su(S^+_{\mathfrak{s}}) \tensor ((\u(E_1) \oplus \u(E_2))
\cap \su(E)))
\ .
\end{equation}
Now decompose the spinor as $\Psi= \Psi_1 + \Psi_2$, where $\Psi_i \in
\Gamma(X,S^+_{\mathfrak{s}} \tensor E_i)$. Recall that the quadratic map $\mu_{0,0}$
is defined to be
$\mu_{0,0}(\Psi) = \mu_{0,0}(\Psi,\Psi)$, where on the right we mean the bilinear map $\mu_{0,0}$.
We get
\begin{equation*}
\begin{split}
  \mu_{0,0}(\Psi,\Psi)
  	= & \mu_{0,0}(\Psi_1,\Psi_1) + \mu_{0,0}(\Psi_1, \Psi_2) + \mu_{0,0}(\Psi_2, \Psi_1) +
\mu_{0,0}(\Psi_2,\Psi_2) \ .
\end{split}
\end{equation*}

By the definition of $\mu_{0,0}$ and by the above
equation (\ref{mu-decomposition}) we see that $\mu_{0,0}(\Psi_1,\Psi_2) =
\mu_{0,0}(\Psi_2,\Psi_1) = 0$. Now from the fact that the bilinear map $\mu_{0,0}$ is
`without zero-divisors' by the above Proposition \ref{no zero divisors} we see that
in each fibre $\Psi_1= 0 $ or $\Psi_2=0$. Suppose we have $\Psi_1(x_0) \neq 0$
for some point $x_0 \in X$. As $\Psi$ is continuous we must have
$\Psi \neq 0$ for all $x$ in some neighbourhood $U$ of $x_0$. Therefore
$\Psi_2 \equiv 0$ on $U$. However, the Dirac equation $D_{\hat{A}} \Psi = 0$
implies that $D_{\hat{A}_1} \Psi_1 = 0$ and that $D_{\hat{A}_2} \Psi_2 = 0$,
where the Dirac operator $D_{\hat{A}_i}: \Gamma(X,S^+_{\mathfrak{s}} \tensor
E_i) \to
\Gamma(X,S^-_{\mathfrak{s}} \tensor E_i)$ is defined to be the composition of
$\nabla_{\hat{A_i},B}: \Gamma(X, S^+_{\mathfrak{s}} \tensor E_i) \to \Omega^1(X,
S^+_{\mathfrak{s}} \tensor
E_i)$ with the Clifford-map $\gamma: T^*X \tensor (S^+_{\mathfrak{s}} \tensor
E_i) \to (S^-_{\mathfrak{s}}
\tensor E_i)$. But for each of these Dirac operators there is a unique
continuation theorem for elements in its kernel by Aronaszajin's Theorem
\cite{A}. Therefore, as $\Psi_2 \equiv 0 $ on $U$, it must vanish identically on
$X$. The general case follows easily by iterating the same argument. \qed

Another important result is the following finiteness property of the fixed-point set inside the moduli space:
%
%

\begin{prop}
The set $\mathscr{B}_{\data}^{S^1}$ has non-empty intersection with the moduli space
$M_{\data}$ at most in the locus of vanishing spinor and in
finitely many components $\mathscr{B}_{\data}^{[F]}$, where $[F]$ runs through
the set of isomorphism classes of proper summands of $E$. 
\end{prop}
{\em Proof:}
We will show that if $[\Psi,\hat{A}] \in M_{\data} \cap
\mathscr{B}_{\data}^{[F]}$,
then $c_1^{\R}(F)$ lies in a bounded set within $H^2(X,\R)$, and 
$\langle c_2(F),[X] \rangle \in \Z$ is bounded also. As $c_1^{\R}(F)$ is in the image of the morphism $H^2(X,\Z) \to
H^2(X,\R)$, it will follow that $c_1(F)$ lies in a finite set.
Also, the morphism $H^4(X,\Z) \to \Z$, given by evaluating a classes $\alpha \in
H^4(X,\Z)$ on the fundamental cycle $[X]$, is an isomorphism. The conclusion is
then that only finitely many pairs $(c_1,c_2) \in H^2(X,\Z) \times H^4(X,\Z)$
can occur as first and second Chern-class of $F$. But on a closed oriented
4-manifold unitary bundles are classified, up to isomorphism, by their first
and second Chern class.

Recall the Chern-Weil formulae for the image of the first and second Chern class
inside $H^*(X,\R) \cong H^*_{dR}(X)$:
\begin{equation}\label{Chern-Weil}
\begin{split}
  c_1^{\R}(E) = & \ \frac{-1}{2\pi i} \left[ \tr F_{\hat{A}} \right] \\
  c_2^{\R}(E) = & \ \frac{-1}{4\pi^2} \left[ \frac{1}{2} \left( \tr{F_{\hat{A}}}
\wedge
\tr{F_{\hat{A}}} - \tr ( F_{\hat{A}} \wedge F_{\hat{A}}) \right) \right] \\
    = & \  \frac{1}{2} \ c_1^{\R}(E)^2 + \frac{1}{8\pi^2} \ [\tr ( F_{\hat{A}}
\wedge
F_{\hat{A}}) ]
\end{split}
\end{equation}
In particular, 
\begin{equation}\label{Chern_2_evaluated}
\begin{split}
\langle c_2(E), [X] \rangle = & \frac{1}{2} \langle c_1(E)^2, [X] \rangle +
\frac{1}{8\pi^2} \int_X \tr (F_{\hat{A}} \wedge F_{\hat{A}}) \\ 
 = & \frac{1}{2} \langle c_1(E)^2, [X] \rangle 
 	+ \frac{1}{8\pi^2} \left( \norm{F_{\hat{A}}^-}_{L^2(X)}^2 
		- \norm{F_{\hat{A}}^+}_{L^2(X)}^2 \right) \ .
\end{split}
\end{equation} 
Let us denote by $\mathscr{H}^2(X,g)$ the vector space of harmonic 2-forms on
$X$. As a subspace of the vector space $\Omega^2(X)$ it is a real inner product space, and the de Rham isomorphism theorem states that it is canonically isomorphic to the second de Rham cohomology group $H^2(X,\R)$, the isomorphism sending simply a harmonic form to its de Rham cohomolgy class that it defines. This way we get thus an inner product on the cohomology group $H^2(X;\R)$. In particular, for a class $[\omega] \in H^2(X;\R)$ we have 
\[
\norm{[\omega]}^2 \leq \norm{\omega}^2_{L^2(X)} = \int_X \omega \wedge * \omega \ .
\]
It follows that if we want to give a bound on the class $c_1(F) = [\tr
F_{\hat{A}_1}]$ inside $H^2(X;\R)$, $\hata_1$ being a connection on $F$, it will be enough to bound the norm
\[
  \norm{\tr F_{\hat{A}_1}}_{L^2(X)} \ .
\]

Now by the assumption that $[\Psi,\hat{A}] \in M_{\data} \cap
\mathscr{B}_{\data}^{[F]}$ we have a connection $\hat{A}$ on $E$ that reduces to $\hat{A}_1 \oplus
\hat{A}_2$ according to the splitting $E = F \oplus F^{\perp}$. We therefore get the decomposition
\[
  F_{\hat{A}} = \begin{pmatrix} F_{\hat{A}_1} & 0 \\ 0 & F_{\hat{A}_2}
\end{pmatrix} \ .
\]
In particular,
\begin{equation}\label{bound_FA1_FA}
  \abs{F_{\hat{A}_1}}^2 \leq \abs{F_{\hat{A}}}^2\ .
\end{equation}
On the other hand, the Lie algebra $\u(F)$ has an orthogonal decomposition
$\u(F)= \su(F) \oplus i \R \id_{F}$. Accordingly we have
\begin{equation}\label{decomp_innerproduct_on_uE}
  \abs{F_{\hat{A}_1}}^2 = \abs{F_{A_1}}^2 + \frac{1}{n} \abs{\tr
F_{\hat{A}_1}}^2 \ ,
\end{equation}
where $n$ equals the rank of $F$.
Therefore 
\begin{equation}\label{bound_trFA1}
\begin{split}
\abs{\tr (F_{\hat{A}_1})}^2 \leq &   n \abs{F_{\hat{A}_1}}^2 \\
  \leq  & n \abs{F_{\hat{A}}}^2 \\
  =  & n \left( \abs{F_A}^2 + \frac{1}{N} \abs{\tr F_{\hat{A}}}^2 \right) \ .
\end{split}
\end{equation}
Recall that $\gamma: \Lambda^2_+ \tensor \su(E) \to \su(S^+_{\mathfrak{s}})
\tensor \su(E)$ is
an isometry, up to universal constant rescaling. The self-dual part
$F_{\hat{A}}^+$ of the trace-free curvature is bounded through the
$PU(N)$-monopole equations (\ref{PUN-equations}) and the a priori bound on the 
spinor (\ref{apriori-bound}):
\begin{equation}\label{FAplusbound}
\begin{split}
\abs{F_A^+}^2 \leq & \ C_0 \ \abs{\mu(\Psi)}^2 \\
  \leq & \  C_0 \ C_1 \ \abs{\Psi}^2 \\
  \leq & \ C_0 \ C_1 \ \frac{4K^2}{c^4} \ .
\end{split}
\end{equation}
Here $C_0$ is a constant that depends on $\gamma$ and $C_1$ is a constant that
depends on the map $\mu$ only. The constant $K$ depends on the metric, our fixed
$Spin^c$-connection $B$ and the fixed connection $\theta$ in the determinant
line bundle of $E$. On the other hand we have $\tr F_{\hat{A}} = F_{\theta}$,
and therfore $\tr F_{\hat{A}}^+$ is bounded also. Therefore we get
\begin{equation} \label{L^2-bound_sd_part}
\begin{split}
\norm{F_{\hat{A}}^+}^2_{L^2(X)} \leq \ C_0 \ C_1 \ \frac{4K^2}{c^4} \
\text{vol}(X) + \norm{F_{\theta}^+}_{L^2(X)}^2 
\end{split}
\end{equation}
Now using the Chern-Weil formula (\ref{Chern_2_evaluated}), we get a bound on
the $L^2$-norm of the anti-self-dual part of $F_{\hat{A}}$ as well:
\begin{equation} \label{L^2-bound_asd_part}
\begin{split}
\norm{F_{\hat{A}}^-}^2_{L^2(X)} \leq & 8 \pi^2 \left( \langle c_2(E), [X]
\rangle 
	- \frac{1}{2} \langle c_1(E)^2, [X] \rangle \right) +
\norm{F_{\hat{A}}^+}^2_{L^2(X)} \\
\end{split}
\end{equation}
The equations (\ref{L^2-bound_sd_part}) and (\ref{L^2-bound_asd_part}) together
imply that the $L^2$-norm of the whole curvature $F_{\hat{A}}$ is bounded:
\begin{equation} \label{L^2-bound_FA}
\begin{split}
\norm{F_{\hat{A}}}^2_{L^2(X)} \leq & 8 \pi^2 \left( \langle c_2(E), [X] \rangle 
	- \langle c_1(E)^2, [X] \rangle \right) \\ & + \ C_0 \ C_1 \ 
\frac{8K^2}{c^2} \ \text{vol}(X) + 2 \norm{F_{\theta}^+}_{L^2(X)}^2 \ .
\end{split}
\end{equation} 
The bound is given by expressions that depend on the metric $g$, the
$Spin^c$ connection $B$ and the connection $\theta$ in the determinant line
bundle, as well as on some constants related to the $\mu$-map and $\gamma$, so
it is a uniform bound on $M_{\data}$. Let us denote by $K_2$ the left
hand
side of the last inequality. From the inequality (\ref{bound_trFA1}) above we
get the desired bound on the $L^2$-norm of $\tr F_{\hat{A_1}}$:
\begin{equation*}
  \norm{\tr F_{\hat{A}_1}}^2_{L^2(X)} \leq n \ K_2 \ .
\end{equation*}
This proves the first assertion, $c_1(F)$ belongs to a finite
subset of the second cohomology group $H^2(X,\Z)$.

The Chern-Weil-formulae for $F$ give a bound on the absolute value of the
second Chern number of $F$:
\begin{equation*}
\begin{split}
 \abs{\langle c_2(F), [X] \rangle } = & \ \frac{1}{2} \abs{ \langle c_1(F)^2,
[X]
\rangle}  - \frac{1}{8 \pi^2} \norm{F_{\hat{A}_1}^+}^2_{L^2(X)} + \frac{1}{8
\pi^2} \norm{F_{\hat{A}_1}^-}^2_{L^2(X)} \\
  \leq & \ \frac{1}{2} \abs{ \langle c_1(F)^2, [X] \rangle}  + \frac{1}{8 \pi^2}
\norm{F_{\hat{A}_1}}^2_{L^2(X)} \ .
 \end{split}
\end{equation*}
As we have $\abs{F_{\hat{A}_1}} \leq \abs{F_{\hat{A}}}$ we get 
\[
  \norm{F_{\hat{A}_1}}^2_{L^2(X)} \leq  \norm{F_{\hat{A}}}^2_{L^2(X)} ,
\]
and the latter is uniformly bounded. As we already have proved that $c_1(F)$
is uniformly bounded in $H^2(X,\R)$ we get that $\abs{ \langle c_1(F)^2, [X]
\rangle}  $ is also.  This implies eventually that $ \abs{\langle c_2(F), [X]
\rangle }$ is uniformly bounded in $\Z$. \qed

\begin{remark}
The $S^1$-action extends naturally to the Uhlenbeck-com\-pacti\-fication
\[
\bar{M}_{\data} \ \subset \ I M_{\data} = \coprod_{k \geq 0}
M_{\data_{-k}} \times Sym^k(X) \ .
\]
As only the second Chern-classes decrease in the Uhlenbeck-com\-pacti\-fication,
the determinant line bundle being preserved, we see by the particular estimates
in the preceeding proof that not only the main-stratum $M_{\data}$
contains only finitely many $S^1$-fixed point components, but the same is even
true for the whole Uhlenbeck-compactification of $M_{\data}$.
\end{remark}
%
In the sequal we shall denote by $M_{\data}^{S^1}$ the intersection
$\mathscr{B}_{\data}^{S^1} \cap M_{\data}$, as well as by
$M_{\data}^{[F]}$ the intersection of $\mathscr{B}_{\data}^{[F]} \cap
M_{\data}$. Also, $M_{\data}^*$ shall denote the
intersection of $\mathscr{B}_{\data}^*$ with $M_{\data}$, and 
$M^*{}_{\data}^{S^1}$, $M_{\data}^*{}^{[F]}$ the
respective intersections with the fixed point set and the given fixed point set
component.

\subsection{Monopole equations for configurations mapping to the fixed point
set}
Above we have pointed out that for describing the component of the fixed point
set $\mathscr{B}_{\data}^{[F]}$ determined by the isomorphism class of a proper
subbundle $[F]$ of $E$, it is useful to keep a representative $F$ fixed. We did
then describe the component $\mathscr{B}_{\data}^{[F]}$ as the image via
$i_F$ of the space $\mathscr{B}_{\data}^{F \oplus F^\perp}$ which is easier
to handle with.
It will turn out that this way we also get a convenient description of
$M_{\data}^{[F]}$, which we define to be the intersection of $\bonf_{\data}^{[F]}$ with the moduli space $M_{\data}$. 

Let us write down explicitly the
monopole equations which are satisfied by a representative $(\Psi,\hat{A})$
having the property that there is a $\hat{A}$-parallel decomposition of $E$
into $F \oplus F^\perp$, with $\Psi$ a section of $S_{\mathfrak{s}}^+ \tensor
F$, and $\hat{A}$ splitting as $\hat{A}_1 \oplus \hat{A}_2$. Recall that $\det(\hat{A}) = \det(\hata_{1}) \tensor \det(\hata_2)$ is the fixed connection $\theta$ in
the determinant line bundle $\det(E)$. We then have 
\[
  F_{A} = F_{\hata} - \, \frac{1}{N} \, \tr F_{\hata} \, \id_E = \begin{pmatrix} F_{\hat{A}_1} -\frac{1}{N} \, F_\theta & 0 \\ 0 & F_{\hat{A}_2} - \frac{1}{N} \, F_\theta
\end{pmatrix} \   ,
\]
according to the splitting $E = F \oplus F^\perp$, and also

\begin{equation*}
\begin{split}
\mu_{0,0}(\Psi) & = \begin{pmatrix} \mu^{F}_{0,1}(\Psi) - \ \frac{1}{N} \ \tr \ 
\mu_{0,1}^{F}(\Psi) \ \id_{F} & 0 \\ 
0 & -\frac{1}{N} \ \tr \ \mu_{0,1}^{F}(\Psi) \ \id_{F^\perp} \end{pmatrix} \\
	& = \begin{pmatrix} \mu^{F}_{0,1-\frac{n}{N}}(\Psi) & 0 \\ 
0 & -\frac{1}{N} \ \tr \ \mu_{0,1}^{F}(\Psi) \ \id_{F^\perp} \end{pmatrix} \ . \end{split}  \end{equation*}
The $PU(N)$-monopole equations \ref{PUN-equations} for the pair $(\Psi, \hat{A}_1 \oplus \hata_2)$ then read

\begin{equation}\label{split_monopoles prov}
\begin{split}
 D_{\hat{A}_1} \Psi & =  0 \\
 \gamma\left(F_{\hat{A}_1}^+ \right) -
\mu_{0,1-\frac{n}{N}}^F(\Psi) & = \frac{1}{N} \, \gamma \klammer{ F_\theta^+ \ \id_{F}} \\
 \gamma\left(F_{\hat{A}_2}^+ \right) +
\frac{1}{N} \ \tr \ \mu_{0,1}^F(\Psi) \  \id_{F^\perp} & =  \frac{1}{N}\, \gamma
\klammer{ F_\theta^+ \ \id_{F^\perp}} \ . \\
\end{split}
\end{equation}
Here the terms in the second equation are sections of the bundle $\su(S_\mathfrak{s}^+)\tensor_\R \u(F)$ and the terms in the third equation are sections of $\su(S_\mathfrak{s}^+) \tensor_\R \u(F^\perp)$. 
There is Lie algebra decompositions $\u(F) = \su(F) \oplus i \R$ and correspondingly for $\u(F^\perp)$. It turns out that the `$i \R$' component of the second and the third equation are equivalent. Indeed, taking the trace (with respect to the factor $\u(F)$ in 
$\su(S_\mathfrak{s}^+)\tensor_\R \u(F)$, and correspondingly for $\u(F^\perp)$) of the second and the third equation, and using the fact that 
\[
\tr(F_{\hata_1}) + \tr(F_{\hata_2}) = F_\theta \ ,
\]
this follows from a simple computation. Therefore the system of equations (\ref{split_monopoles prov}) above is equivalent to the same system where we take as the third equation only the component of $\su(F^\perp)$ according to $\u(F^\perp) = \su(F^\perp) \oplus i \R $. Thus the $PU(N)$ monopole equations are therefore equivalent to 
\begin{equation}\label{split monopoles}
\begin{split}
 D_{\hat{A}_1} \Psi & =  0 \\
 \gamma\left(F_{\hat{A}_1}^+ \right) -
\mu_{0,1-\frac{n}{N}}^F(\Psi) & = \frac{1}{N} \, \gamma \klammer{ F_\theta^+ \ \id_{F}} \\
 \gamma\left(F_{A_2}^+ \right)  & =  0 \ . \\
\end{split}
\end{equation}
We summarise this computation in the following:
\begin{prop}\label{decoupled equations}
Suppose the configuration $(\Psi,\hata) \in \conf_{\data}$ has reducible connection part $\hata = \hata_1 \oplus \hata_2$ according to $E = F \oplus F^\perp$, and that the spinor part $\Psi$ is a section of $S^+ \tensor F$ (compare proposition \ref{one_component_vanishes}). Then the $PU(N)$ monopole equations for $(\Psi,\hata)$ are equivalent to the system (\ref{split monopoles}).  
\end{prop}


\begin{definition}
We shall denote by $M_{\data}^{F \oplus F^\perp} \subseteq
\mathscr{B}_{\data}^{F \oplus F^\perp}$ the moduli space space of solutions
$(\Psi,\hat{A}_1,\hat{A}_2) $ to the above equations
(\ref{split monopoles}) modulo the gauge group $\G_{F \oplus F^\perp}$.
As usually, we denote by $M_{\data}^{*\ F \oplus F^\perp}$ the subspace
of those elements whose representatives have zero-dimensional stabiliser.
\end{definition}

\begin{prop}\label{iF_modulispace}
The map $i_F:\mathscr{B}_{\data}^{F\oplus F^\perp} \to
\mathscr{B}_{\data}^{[F]}$ maps the moduli space $M_{\data}^{F \oplus
F^\perp}$ onto the fixed point set component $M_{\data}^{[F]}$ inside
the moduli space. It maps the set $M_{\data}^*{}^{F \oplus F^\perp}$
bijectively onto $M^*{}_{\data}^{[F]}$.

\end{prop}
{\em Proof:} The fact that the map is onto is an immediate consequence of the above Proposition \ref{decoupled equations} and the definition of $M_{\data}^{[F]}$. 
For the remaining claim we will show that we can apply the above Proposition \ref{Prop_of_iF}. First, we shall observe that
 if $[\Psi,\hat{A}_1,\hat{A}_2]$ belongs to $M_{\data}^{*\ F \oplus
F^\perp}$, then the connections $\hat{A}_1$ and $\hat{A}_2$ are indeed
irreducible. Obviously $\hata_2$ has to be irreducible, but suppose $\hat{A}_1$ were reducible. We would have a $\hat{A}_1$ - parallel orthogonal decomposition $F = F_1 \oplus
F_2$, with $\hat{A}_1$ splitting accordingly, $\hat{A}_1 = \hat{A}_{11} \oplus
\hat{A}_{12}$. 

Let us write $\Psi=\Psi_1 + \Psi_2 \in
\Gamma(X,S^+_{\mathfrak{s}} \tensor (F_1 \oplus F_2))$ for the corresponding decomposition of
the spinor. We claim that either $\Psi_1 = 0$ or $\Psi_2 = 0$.
In fact, $(\Psi,\hata_1)$ solves the first two of the equations (\ref{split monopoles}). The map $\mu_{0,\tau}^F$ is `without zero-divisors' by the above Proposition \ref{no zero divisors}. With this fact the conclusion follows exactly like in the proof of Proposition \ref{one_component_vanishes}. But then the configuration
$(\Psi,\hat{A}_{11}\oplus \hat{A}_{12},\hat{A}_2)$ must have
positive-dimensional stabiliser inside $\G_{F \oplus F^\perp}$, and the element 
$[\Psi,\hat{A}_1,\hat{A}_2]$ would not belong to $M_{\data}^{*\ F
\oplus
F^\perp}$. Therefore $[\Psi,\hat{A}_1,\hat{A}_2]$ belongs to the set
$\mathscr{B}_{\data}^{*,irr \ F \oplus
F^\perp}$ and we can apply Proposition $\ref{Prop_of_iF}$ for getting
injectivity. Furthermore it is easy to see that the parametrisation $i_F$
maps $M_{\data}^*{}^{F \oplus F^\perp}$
onto $M^*{}_{\data}^{[F]}$.
\qed

The first two equations of (\ref{split monopoles}) are $U(n)$- monopole equations for $(\Psi,\hata_1)$ as discussed in \cite{Z_un} and the third equation is the anit-self-duality equation for the $PU(N-n)$ connection $A_2$. We recall the definition of the corresponding moduli spaces. In the terminology of \cite{Z_un} the space of configurations $(\Psi,\hata_1) \in \conf_{\dataf}^U$ which satisfy the first two equations of (\ref{split monopoles}) quotiented by the gauge-group $\mathscr{G}_F$ of unitary automorphisms of $F$ is denoted by $M_{\dataf}(1-\frac{n}{N},\frac{1}{N}F_\theta^+)$, and called the moduli space of $U(n)$ - monopoles of parameter $1-\frac{n}{N}$, and perturbed by the self-dual two-form $\frac{n}{N} F_\theta^+$. 

We shall further denote by $M_{F}^{asd}$ the moduli space of anti-self-dual $PU(n)$ - connections in $F$ which is defined to be the space of $PU(n)$- connections $A \in \mathscr{A}^{PU}(F)$ in $F$ which satisfy the equations $F_A^+ = 0$, quotiented by the action of the gauge-group $\G_{F}$ of special unitary automorphisms of $F$. 

In Proposition \ref{decoupled equations} we observed that the monopole equations for the configuration $(\Psi,\hata_1 \oplus \hata_2 \in \conf_{\data}^{F \oplus F^\perp})$ were decoupling into a $U(n)$-monopole equation for $(\Psi,\hata_1)$ and an anti-self-duality equation for the $PU(N-n)$- connection $A_2$. Given the product description $\bonf_{\data}^{F \oplus F^\perp} \cong \bonf_{\dataf}^U \times \bonf_{F^\perp}^{PU}$ in Proposition \ref{map_h} it is now not surprising that there is a corresponding product description for the moduli spaces, as we shall see in the next Proposition. For better notation we shall write $M_{\dataf}^U$ for the moduli space $M_{\mathfrak{s},F}\left(1-\frac{n}{N},\frac{1}{N} \, F_\theta^+\right)$. 
\begin{prop}\label{moduli_space_product}
  Restricting the bijection $h$ of Proposition \ref{map_h} above to the moduli
space $M_{\data}^{F \oplus F^\perp}$ we get an induced bijection
\begin{equation*}
h|_{M} : M_{\data}^{F \oplus F^\perp} \stackrel{\cong}{\to}
\
M_{\mathfrak{s},F}^U \times M^{asd}_{F^\perp} 
\end{equation*}
Together with the map $i_F$ we thus get the parametrisation of the fixed point
set component $M_{\data}^{[F]}$ of the moduli space as the
product of
a moduli space of $U(n)$-monopoles with the moduli space of
$ASD-PU(N-n)$-connections. In particular, for the irreducible parts we get a
bijection
\begin{equation*}
i_F \circ h|_M^{-1} : M_{\mathfrak{s},F}^{* \ U} \times
M^{* \ ASD}_{\ F^\perp}  \stackrel{\cong}{\to}
M_{\data}^{* \ [F]} \ .
\end{equation*}
\end{prop}
This follows from Proposition \ref{iF_modulispace} and Proposition \ref{map_h},
where it is easily checked that $h|_M^{-1}$ maps the `irreducibles'
$M_{\dataf}^{* \, U} \times M_{F^\perp}^{* \, ASD}$ onto the
corresponding `irreducibles' $M_{\data}^{* \ F\oplus F^\perp}$. 
\qed 

The whole discussion is now summarised in 

\begin{theorem}\label{S1fixedpointset_main_text}
The fixed point set under the above circle-action $r$ on the moduli space $M_{\data,\theta}$ of $PU(N)$  monopoles is given as the union of the moduli space $M_E^{asd}$ of
anti-self-dual $PU(N)$ connections in $E$ and a finite union  
\[
  \bigcup_{[F] \subseteq E} M_{\data}^{[F]} 
\]
 of components $M_{\data}^{[F]}$ indexed by a finite number of isomorphism classes $[F]$  of
proper subbundles of $E$. The spaces $M_{\data}^{[F]}$ are given as
follows: An element $[\Psi,\hat{A}]$ belongs to $M_{\data}^{[F]}$ if 
for each representative $F \in [F]$ there is a
representative $(\Psi,\hat{A}) \in [\Psi,\hat{A}]$ such that $F$ is a
$\hat{A}$-invariant proper subbundle of $E$ and the spinor $\Psi$ is a section of the
proper subbundle $S^+_\mathfrak{s} \tensor F$ of $W^+_{\data}=S^+_\mathfrak{s} \tensor E$.

 Furthermore, if $X$ is simply connected, then there is a parametrisation of
this
space $M_{\data}^{[F]}$ as a product
\[
  M_{\mathfrak{s},F}^U \times M_{F^\perp}^{asd} \to
M_{\data}^{[F]} \ .
\]
Here $M_{\mathfrak{s},F}^{U}$ is a moduli space of
$U(n)$ monopoles, with $n$ being the rank of $F$, and
$M_{F^\perp}^{asd}$ is the moduli space of $PU(N-n)$ anti-self-dual
 connections in $F^\perp$. This map is a surjection and is a bijection between the open and dense subsets of elements with zero-dimensional stabiliser in the corresponding moduli spaces, 
\[
  M_{\mathfrak{s},F}^{*\ U} \times M_{F^\perp}^{*\ asd} \stackrel{\cong}{\to}
M_{\data}^{*\ [F]} \ .
\]
\end{theorem}

We didn't carry through the topological part of the `cobordism program' here. However, there is evidence that there will be no contributions to the higher rank instanton invariant coming from moduli spaces $M^{[F]}_{\data} \cong M_{\mathfrak{s},F}^U \times M^{asd}_{F^\perp}$ with $\rk(F) > 1$ because the main results of \cite{Z_un} are as follows: First, the moduli spaces of $U(n)$ Seiberg-Witten monopoles with $n > 1$ are empty on K\"ahler surfaces with $b^+_2(X) > 1$ if we perturb with a non-vanishing holomorphic two-form. Second, the $U(n)$ monopole equations can be perturbed in a way so that the moduli spaces of the perturbed equations are empty after perturbing with a generic self-dual two-form if $b^+_2(X) > 0$.

\end{document}